\documentclass[graybox]{svmult}

\usepackage{mathptmx}       
\usepackage{helvet}         
\usepackage{courier}        
\usepackage{type1cm}        
\usepackage{graphicx}        
\usepackage[bottom]{footmisc}

\usepackage{amsmath,amssymb}
\usepackage{tikz}

                       \newcommand{\Z}{\mathbb{Z}}%

\numberwithin{theorem}{section}

\newtheorem{lem}[theorem]{Lemma}
\newtheorem{cor}[theorem]{Corollary}

\newcommand{\R}{\mathbb{R}}%
\newcommand{\C}{\mathbb{C}}%
\newcommand{\T}{\mathbb{T}}
\newcommand{\TT}{\mathcal{T}}
\newcommand{\OO}{\mathcal{O}}

\newcommand{\clsp}{\overline{\operatorname{span}}}

\newcommand{\Aut}{\operatorname{Aut}}

\newcommand\id{\mathrm{id}}

\begin{document}

\title*{Equilibrium states on graph algebras}
\author{Astrid an Huef and Iain Raeburn}
\institute{Astrid an Huef and Iain Raeburn \at Department of Mathematics and Statistics, University of Otago, PO Box 56, Dunedin 9054, New Zealand. \at \email{\{astrid, iraeburn\}@maths.otago.ac.nz}
}
%
%
\maketitle

{\abstract{We consider operator-algebraic dynamical systems given by actions of the real line on unital $C^*$-algebras, and especially the equilibrium states (or KMS states) of such systems. We are particularly interested in systems built from the gauge action on the Toeplitz algebra and graph algebra of a finite directed graph, and we describe a complete classification of the KMS states obtained in joint work with Laca and Sims. We then discuss applications of these results to Cuntz-Pimsner algebras associated to local homeomorphisms, obtained in collaboration with Afsar. Thomsen has given  bounds on the range of inverse temperatures at which KMS states may exist. We show that Thomsen's bounds are sharp.}}

\section{Introduction}

There has recently been a renewal of interest in the equilibrium states (the KMS states) of operator-algebraic dynamical systems consisting of an action of the real line  (the \emph{dynamics}) on a $C^*$-algebra. There has been particular interest in systems involving graph algebras and their Toeplitz extensions \cite{EL, KW, T, CL}. 

Very satisfactory results have been obtained for sytems associated to finite directed graphs, and we now have concrete descriptions of the simplices of KMS$_\beta$ states on the Toeplitz algebras at all inverse temperatures $\beta$ \cite{aHLRS1, aHLRS2}. Here we review these results, and discuss some surprising applications to work of Thomsen on systems involving the Cuntz-Pimsner algebras of local homeomeorphisms \cite{Th1}. One main conclusion of our recent work with Afsar \cite{AaHR} is that lower and upper bounds for the possible inverse temperatures given by Thomsen are sharp. For these applications we do not need the full strength of the general results in \cite{aHLRS2}, and in this article we describe a more direct approach.

\section{The Toeplitz algebra of a graph}

We suppose that $E=(E^0,E^1,r,s)$ is a finite directed graph. A \emph{Toeplitz-Cuntz-Krieger $E$-family} consists of mutually orthogonal projections $\{P_v:v\in E^0\}$ and partial isometries $\{S_e:e\in E^1\}$ in a $C^*$-algebra such that $S_e^*S_e=P_{s(e)}$ for every $e\in E^1$ and 
\begin{equation}\label{TCK}
P_v\geq \sum_{r(e)=v}S_eS_e^*\text{ for every $v\in E^0$}
\end{equation} 
(where we interpret an empty sums as $0$). Since the vertex projections are mutually orthogonal, the relation \eqref{TCK} implies that the range projections $S_eS_e^*$ are also mutually orthogonal. (See \cite[Corollary~1.2]{aHLRS1}, for example.)

For $n\geq 2$, we write 
\[
E^n:=\big\{\mu=\mu_1\mu_2\cdots\mu_n:s(\mu_i)=r(\mu_{i+1})\text{ for $1\leq i<|\mu|:=n$}\big\}
\]
for the set of paths of length $n$ in $E$, and note that $S_\mu:=S_{\mu_1}S_{\mu_2}\cdots S_{\mu_n}$ is a partial isometry for every $\mu\in E^n$. We write $E^*:=\bigcup_{n\geq 0} E^n$ for the set of finite paths. Then for $\mu,\nu,\alpha,\beta\in E^*$ we have the product formula
\begin{equation}\label{prodform}
(S_\mu S_\nu^*)(S_\alpha S_\beta^*)=\begin{cases}
S_{\mu\alpha'}S_\beta^* & \text{if $\alpha=\nu\alpha'$}\\
S_\mu S_{\beta\nu'}^* &\text{if $\nu=\alpha\nu'$}\\
0&\text{otherwise.}\end{cases}
\end{equation}

The \emph{Toeplitz algebra} $\TT C^*(E)$ of $E$ is the $C^*$-algebra generated by a universal Toeplitz-Cuntz-Krieger $E$-family $(p,s)$. The product formula \eqref{prodform} implies that the elements $\{s_\mu s_\nu^*:\mu,\nu\in E^*\}$ span a $*$-subalgebra of $\TT C^*(E)$, and hence we have
\[
\TT C^*(E)=\clsp\{s_\mu s_\nu^*:\mu,\nu\in E^*\}.
\]
The quotient of $\TT C^*(E)$ by the ideal generated by the gap projections
\[
\Big\{p_v-\sum_{r(e)=v}s_es_e^*:v\in E^0\Big\}
\]
is the usual graph algebra or Cuntz-Krieger algebra $C^*(E)$.

For every graph $E$ there is a canonical Toeplitz-Cuntz-Krieger $E$-family $(Q,T)$ on the \emph{finite-path space} $\ell^2(E^*)$, characterised by the following actions on the usual orthonormal basis $\{h_\mu:\mu\in E^*\}$: 
\[
Q_v h_\mu=\begin{cases} h_{\mu}&\text{if $v=r(\mu)$}\\
0&\text{otherwise} \end{cases}
\quad\text{and}\quad T_e h_\mu=\begin{cases} h_{e\mu}&\text{if $s(e)=r(\mu)$}\\
0&\text{otherwise.} \end{cases}
\]
The universal property of $\TT C^*(E)$ then gives a representation $\pi_T=\pi_{Q,T}$ of $\TT C^*(E)$ on $\ell^2(E^*)$ such that $\pi_T(p_v)=Q_v$ and $\pi_T(s_e)=T_e$; we call $\pi_T$ the \emph{finite-path representation}. The gap projections $Q_v-\sum_{r(e)=v}T_eT_e^*$ are the projections on $\C h_v$, and hence are all nonzero. Thus the uniqueness theorem for Toeplitz algebras \cite[Corollary~4.2]{FR} implies that $\pi_T$ is faithful.

There is a gauge action $\gamma:\T\to \Aut \TT C^*(E)$ such that $\gamma_z(p_v)=p_v$ and $\gamma_z(s_e)=zs_e$, and this induces the usual gauge action on the quotient $C^*(E)$. We are interested in the dynamics $\alpha:\R\to \Aut \TT C^*(E)$ given by $\alpha_t=\gamma_{e^{it}}$, and its analogue on $C^*(E)$. In particular, we wish to study the KMS states for this dynamics.

\section{KMS states on the Toeplitz algebra}
\label{sec:1}

The spanning elements $s_\mu s_\nu^*$ for $\TT C^*(E)$ are are all analytic for the action $\alpha$. Hence if $\phi$ is a KMS$_\beta$ state on $(\TT C^*(E),\alpha)$, we have
\begin{align*}
\phi(s_\mu s_\nu^*)&=\phi(s_\nu^*\alpha_{i\beta}(s_\mu))=e^{-\beta|\mu|}\phi(s_\nu^* s_\mu)\\
&=e^{-\beta(|\mu|-|\nu|)}\phi(s_\mu s_\nu^*).
\end{align*}
So $\phi(s_\mu s_\nu^*)\not=0\Longrightarrow |\mu|=|\nu|$, and then 
\[
\phi(s_\mu s_\nu^*)\not=0\Longrightarrow\big(s_\nu^* s_\mu\not=0\text{ and }|\mu|=|\nu|\big)\Longrightarrow\mu=\nu.
\]

Now a routine computation using the product formula \eqref{prodform} gives the following:

\begin{lem}\label{nascKMS}\cite[Proposition~2.1]{aHLRS1} A state $\phi$ on $\TT C^*(E)$ is KMS$_\beta$ for $\alpha$ if and only if 
\[
\phi(s_\mu s_\nu^*)= 
\delta_{\mu,\nu}e^{-\beta|\mu|}\phi(p_{s(\mu)}).
\]
\end{lem}

Suppose $\phi$ is a KMS$_\beta$ state on $(\TT C^*(E),\alpha)$. Then for each $v\in E^0$ the Toeplitz-Cuntz-Krieger relation gives
\begin{align}\label{TCKineq}
\phi(p_v)\geq\sum_{r(e)=v}\phi(s_es_e^*)=\sum_{r(e)=v}e^{-\beta}\phi(p_{s(e)}),
\end{align}
where we interpret the empty sum as $0$ if $v$ is a source. The \emph{vertex matrix} of $E$ is the $E^0\times E^0$ integer matrix $A$ with entries \[
A(v,w)=|r^{-1}(v)\cap s^{-1}(w)|.\]
We can rewrite the inequality \eqref{TCKineq} as
\begin{equation}\label{subinv}
e^{\beta}\phi(p_v)\geq\sum_{w\in E^0}\;\sum_{r(e)=v,\;s(e)=w}\phi(p_w)=\sum_{w\in E^0}A(v,w)\phi(p_{w}),
\end{equation}
so $m=(m_v):=(\phi(p_v))\in [0,\infty)^{E^0}$ satisfies $Am\leq e^{\beta}m$; we say that $m$ is a \emph{subinvariant vector} for $A$. 

If $\phi$ factors through a KMS$_\beta$ state of $C^*(E)$, then we have equality throughout \eqref{TCKineq} and \eqref{subinv}, and $m$ satisfies $Am= e^{\beta}m$. If $E$ is strongly connected, $A$ is irreducible and $e^\beta$ has to be the Perron-Frobenius eigenvalue of $A$. The Perron-Frobenius theorem then says many things: thus in particular we know that the eigenvalue is the spectral radius $\rho(A)$, that the eigenspace is one-dimensional, and that there is an eigenvector with positive entries (see \cite[Theorem~2.6]{DZ} or \cite[Theorem~1.6]{seneta}). Since $\phi$ is a state, we have
\[
1=\phi(1)=\sum_v \phi(p_v)=\sum_v m_v, 
\]
and hence $m=\big(\phi(p_v)\big)_{v\in E^0}$ is the unique eigenvector in $(0,\infty)^{E^0}$ with $\|m\|_1=1$. The formula in Lemma~\ref{nascKMS}  says that $\phi(s_\mu s_\nu^*)= 
\delta_{\mu,\nu}e^{-\beta|\mu|}m_{s(\mu)}$ for all $\mu,\nu\in E^*$, so the vector $m$ completely determines the state $\phi$. Thus we recover the following elegant result of Enomoto, Fujii and Watatani \cite{EFW}:

\begin{theorem} Suppose that $E$ is a strongly connected graph with vertex matrix $A$. Then $(C^*(E),\alpha)$ has at most one KMS state. This state has inverse temperature $\ln\rho(A)$, where $\rho(A)$ is the spectral radius of $A$.
\end{theorem}

They proved existence of the KMS$_{\ln\rho(A)}$ state too, but we'll get to that.

For states on $\TT C^*(E)$, the vector $m$ only satisfies the subinvariance relation $Am\leq e^{\beta}m$, but when $A$ is irreducible Perron-Frobenius theory has things to say about this too. For example: 
\begin{itemize}
\item If $Am\leq e^{\beta}m \text{ and }\beta=\ln\rho(A)$, then  $Am= e^{\beta}m$, so that $m$ is the Perron-Frobenius eigenvector.
\item Suppose that $Am\leq e^{\beta}m$. Then $Am\not=e^{\beta}m\Longleftrightarrow\beta>\ln\rho(A)$.
\end{itemize}
This suggests that we look more carefully at $\beta$ larger than the \emph{critical inverse temperature} $\beta_c:=\ln\rho(A)$. 

So we consider $\beta>\ln\rho(A)$. We find it interesting that, although we were motivated to do so by the Perron-Frobenius theory, which applies only when $E$ is strongly connected, the following analysis does not require any connectivity hypothesis on $E$. Thus we consider an arbitrary finite directed graph $E$, which could have sinks or sources, and a KMS$_\beta$ state $\phi$ on $\TT C^*(E)$. 

Take $m=\big(\phi(p_v)\big)$ as before. Then  $\epsilon:=(1-e^{-\beta}A)m$ has nonnegative entries, not all $0$. Since $e^{\beta}>\rho(A)$, $e^{\beta}$ is not in the spectrum of $A$, and $1-e^{-\beta}A$ is invertible. Thus we can recover $m$ as $(1-e^{-\beta}A)^{-1}\epsilon$. Our main  point is that we can describe geometrically the set of $\epsilon\in [0,1]^{E^0}$ which arise from unit vectors $m$ in $\ell^1(E^0)$. For $v\in E^0$, define $y^\beta\in [1,\infty)^{E^0}$ by 
\begin{equation}\label{ybeta}
y^\beta_v:=\sum_{n=0}^\infty\;\sum_{w\in E^0} e^{-\beta n}A^n(w,v)=(1-e^{-\beta}A)^{-1}\delta_v;
\end{equation}
the series converges because $\sum_ne^{-\beta n}A^n$ converges in the operator norm with sum $(1-e^{-\beta}A)^{-1}$. Then $m:=(1-e^{-\beta}A)^{-1}\epsilon$ has $\|m\|_1=1$ if and only if 
\[
1=\epsilon\cdot y^\beta:=\sum_{v\in E^0}\epsilon_vy^\beta_v
\]
(see Theorem~3.1(a) of \cite{aHLRS1}).

Then the main theorem of \cite{aHLRS1} says:

\begin{theorem}\label{paraKMS}
Suppose $E$ is a finite graph with vertex matrix $A$, and $\beta>\ln \rho(A)$. Suppose $\epsilon\cdot y^\beta=1$. Then there is a KMS$_\beta$ state $\phi_\epsilon$ of $\TT C^*(E)$ such that 
\[
\phi_\epsilon(p_v)=\big((1-e^{-\beta}A)^{-1}\epsilon\big)_v\quad\text{for all $v\in E^0$.}
\]
The map $\epsilon\mapsto \phi_\epsilon$ is an isomorphism of $\Delta_\beta=\{\epsilon:\epsilon\cdot y^\beta=1\}$ onto the simplex of KMS$_\beta$ states of $(\TT C^*(E),\alpha)$.
\end{theorem}

We proved existence of the KMS state $\phi_\epsilon$ in \cite[Theorem~3.1(b)]{aHLRS1} by a spatial argument using the finite-path representation $\pi_T$ of $\TT C^*(E)$ on $\ell^2(E^*)$. Then surjectivity of $\epsilon\mapsto \phi_\epsilon$ amounts to our earlier observation that the subinvariant vector $m=\big(\phi(p_v)\big)_{v\in E^0}$ determines a KMS state $\phi$.

The set  $\Delta_\beta=\{\epsilon:\epsilon\cdot y^\beta=1\}$ parametrising the KMS$_\beta$ states is a simplex in the positive cone $[0,\infty)^{E^0}$ of $\R^{E^0}$ with extreme points on the coordinate axes, and the vector $y^\beta$ is normal to this simplex. As $\beta$ decreases to the critical value $\beta_c=\ln\rho(A)$, the terms in the series on the right-hand side of \eqref{ybeta} get larger, and the simplex contracts towards the origin. 

The preceding analysis does not apply when $\beta=\beta_c=\ln\rho(A)$ is critical, because then the matrix $1-e^{-\beta}A$ need not be invertible. However, we can take a sequence $\beta_n$ decreasing to $\ln\rho(A)$, and use weak* compactness of the state space to get a KMS$_{\ln\rho(A)}$ state of $\TT C^*(E)$ \cite[Proposition~4.1]{aHLRS1}. When $E$ is strongly connected, this is the only KMS$_{\ln\rho(A)}$ state, and we can deduce from Perron-Frobenius that it factors through the graph algebra $C^*(E)$. In particular, we recover the existence of the KMS$_{\ln \rho(A)}$ state, first established by other methods in \cite{EFW}.

We can sum up our discussion as follows:

\begin{cor}\label{subinv2}
Suppose that $E$ is a directed graph with vertex matrix $A$, and that $\beta\in (0,\infty)$ satisfies $\beta\geq \ln\rho(A)$. Then the map $\phi\mapsto m^\phi:=\big(\phi(p_v)\big)$ is a bijection of the set of KMS$_\beta$ states of $(\TT C^*(E),\alpha)$ onto the unit vectors $m$ in $[0,\infty)^{E^0}\subset\ell^1(E^0)$ satisfying the subinvariance relation $Am\leq e^{\beta}m$.
\end{cor}

For $\beta> \ln\rho(A)$, Theorem~\ref{paraKMS} is stronger, because it describes the solutions of the subinvariance relation. But for some applications, such as those in \S\ref{sec:dumb}, we can deal directly with the subinvariance relation in an \emph{ad hoc} manner.

\section{Dumbbell graphs}\label{sec:dumb}

We say that a graph $E$ is \emph{reducible} if it is not strongly connected, or equivalently if its vertex matrix $A$ is not irreducible. For $v,w\in E^0$, we write $v\leq w$ to mean that
\[
vE^*w:=\{\mu\in E^*: r(\mu)=v\text{ and }s(\mu)=w\}
\]
is nonempty (or in other words, that there is a path from $w$ to $v$). Then we define a relation $\sim$ on $E^0$ by 
\[
v\sim w\Longleftrightarrow v\leq w\text{ and }w\leq v.
\]
This is an equivalence relation (it is reflexive because $v\in vE^*v$), and we write $E^0/\!\!\sim$ for the set of equivalence classes.

For each $C\in E^0/\!\!\sim$, we define $A_C$ to be the $C\times C$ matrix obtained by deleting all rows and columns involving vertices not in $C$. We can view $A_C$ as the vertex matrix of the subgraph $E_C:=(C, E^1\cap r^{-1}(C)\cap s^{-1}(C),r,s)$. Each $A_C$ is either a $1\times 1$ zero matrix (if $C$ is a singleton set $\{v\}$ and  there is no loop at $v$) or an irreducible matrix (in which case we call $E_C$ a \emph{strongly connected component} of $E$). It is possible to order the set $E^0$ so that the vertex matrix $A$ is block upper-triangular with diagonal blocks $A_C$ (see \cite[\S2.3]{aHLRS2}), and it follows that $\rho(A)=\max_C\rho(A_C)$.  

Now we consider the KMS states on $\TT C^*(E)$ when $E$ is reducible. There are three situations that we have to deal with:
\begin{itemize}
\item For $\beta>\ln \rho(A)$, Theorem~\ref{paraKMS} applies, and we have a $(|E^0|-1)$-dimensional simplex of KMS$_\beta$ states on $\TT C^*(E)$.
\item For $\beta=\ln\rho(A)$,  we focus on the \emph{critical components} $C\in E^0/\!\!\sim$ that have $\rho(A_C)=\rho(A)$. The relation $\leq $ descends to a well-defined relation on the set of critical components, and then a critical component $C$ is \emph{minimal} if $D$ critical and $D\leq C$ imply $D=C$. The behaviour of the KMS$_{\ln\rho(A)}$ states of $\TT C^*(E)$ depends on the location of the  minimal critical components.
\item Recall that a subset $H$ of $E^0$ is \emph{hereditary} if $v\in H$ and $v\leq w$ imply $w\in H$, and that the Toeplitz algebra $\TT C^*(E\backslash H)$ of the graph 
\[
E\backslash H:=(E^0\backslash H, E^1\cap s^{-1}(E^0\backslash H), r, s)
\]
is a quotient of $\TT C^*(E)$ \cite[Proposition~2.1]{aHLRS2}. For $\beta<\ln\rho(A)$, we consider the hereditary subset $H$ of $E^0$ generated by the critical components (which is also generated by the minimal critical components). If $H$ is not all of $E^0$, then we can apply Theorem~\ref{paraKMS} to $E\backslash H$ and get KMS$_\beta$ states of $\TT C^*(E)$ for $\ln\rho(A_{E\backslash H})<\beta<\ln\rho(A)$ which factor through the quotient map onto $\TT C^*(E\backslash H)$.
\end{itemize}

The first and third situations both require straightforward applications of Theorem~\ref{paraKMS}, and the interesting things happen when $\beta=\ln\rho(A)$ is critical. Then the phrase ``depends on the location of the minimal critical components'' needs clarification. We illustrate its meaning with some examples, which fortunately are enough for the main applications in \cite{AaHR}. The key feature of these examples is that the strongly connected components have just one vertex each. We call such graphs \emph{dumbbell graphs}.

\begin{example} We consider the following graph $E$:
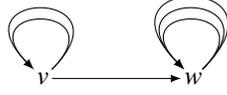
\[
\begin{tikzpicture}
    \node[inner sep=1pt] (v) at (0,0) {$v$};
    \node[inner sep=1pt] (w) at (2,0) {$w$};
    \draw[-latex] (v)--(w);
    \foreach \x in {0,2} {
        \draw[-latex] (v) .. controls +(1.\x,1.\x) and +(-1.\x,1.\x) .. (v);
    }
    \foreach \x in {0,2,4} {
        \draw[-latex] (w) .. controls +(1.\x,1.\x) and +(-1.\x,1.\x) .. (w);
    }
\end{tikzpicture}
\]
for which $\rho(A)=3$. In this example, the hereditary closure of the critical component $\{w\}$ is all of $E^0$, and hence the third situation does not arise.

If we list $E^0=\{w,v\}$, then $A=\big(\begin{smallmatrix}3&1\\0&2\end{smallmatrix}\big)$. For $\beta=\ln \rho(A)=\ln 3$ we have $e^\beta=3$, and the subinvariance relation $Am\leq e^\beta m$ says
\[
Am=\begin{pmatrix}3&1\\0&2\end{pmatrix}\begin{pmatrix}m_w\\m_v\end{pmatrix}=\begin{pmatrix}3m_w+m_v\\3m_v\end{pmatrix}\leq 3\begin{pmatrix}m_w\\m_v\end{pmatrix}.
\]
The only unit vector in $[0,\infty)^{E^0}\subset \ell^1(E^0)$ which satisfies this relation is $m=(1,0)$. Thus Corollary~\ref{subinv2} says there is a unique KMS$_{\ln 3}$ state on $\TT C^*(E)$. This state factors through $C^*(E)$. 
\end{example}

\begin{example}\label{2nddumb}
Next we switch the horizontal arrow, so $E$ is
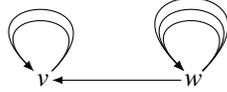
\[
\begin{tikzpicture}
    \node[inner sep=1pt] (v) at (0,0) {$v$};
    \node[inner sep=1pt] (w) at (2,0) {$w$};
    \draw[-latex] (w)--(v);
    \foreach \x in {0,2} {
        \draw[-latex] (v) .. controls +(1.\x,1.\x) and +(-1.\x,1.\x) .. (v);
    }
    \foreach \x in {0,2,4} {
        \draw[-latex] (w) .. controls +(1.\x,1.\x) and +(-1.\x,1.\x) .. (w);
    }
\end{tikzpicture}
\]
With $E^0=\{v,w\}$, $A=\big(\begin{smallmatrix}2&1\\0&3\end{smallmatrix}\big)$, and subinvariance for $\beta=\ln 3$ reduces to $m_v\leq m_w$. This graph also has just one critical component $\{w\}$, but this time $\{w\}$ is hereditary, and the graph $E\backslash H$ has vertex set $\{v\}$, so the third situation kicks in.

We find:
\begin{itemize}
\item For $\beta>\ln 3$, we have a one-dimensional simplex of KMS$_\beta$ states on $\TT C^*(E)$, none of which factor through $C^*(E)$.
\item The simplex of KMS$_{\ln 3}$ states on $\TT C^*(E)$ has extreme points $\phi_v$ ($m_v=m_w=\frac{1}{2}$) and $\phi_w$ ($m_v=0$); only $\phi_w$ factors through a state of $C^*(E)$.
\item For $\ln 2\leq \beta<\ln 3$, there is a unique KMS$_\beta$ state $\phi_v$ on $\TT C^*(E)$ which factors through the quotient map corresponding to the hereditary set $\{w\}\subset E^0$.
\item For $\beta=\ln 2$, the state $\phi_v$ factors through $C^*(E)$.
\item For $\beta<\ln 2$, there are no KMS$_\beta$ states on $\TT C^*(E)$.
\end{itemize}
\end{example}

When a minimal strongly connected component $E_C$ has more than one vertex, we organise the block form for $A$ in three pieces: we take the hereditary closure $H$ of the critical components, and decompose $E^0=(E^0\backslash H)\cup C\cup (H\backslash C)$. The Perron-Frobenius eigenvector for $A_C$ gives a KMS$_{\ln \rho(A)}$ state $\psi_C$ that has $\phi(p_v)=0$ for $v\in H\backslash C$, but has $\phi(p_v)\not=0$ for vertices $v$ such that $vE^*C\not=\emptyset$: the precise formula is given in \cite[Theorem~4.3(a)]{aHLRS2}. Since $\rho(A_{E\backslash H})<\rho(A)$, we can also use Theorem~\ref{paraKMS} to find KMS$_{\ln\rho(A)}$ states on $\TT C^*(E\backslash H)$, and lift them to KMS$_{\ln\rho(A)}$ states of $\TT C^*(E)$. Again the formulas and the complete classification are given in \cite[Theorem~4.3]{aHLRS2}. 

To construct KMS states on the usual graph algebra $C^*(E)$, we need to know which states on $\TT C^*(E)$ factor through $C^*(E)$. Here we hit another subtlety: distinct hereditary sets give distinct ideals in $\TT C^*(E)$ but not necessarily in $C^*(E)$, where the  ideal in $C^*(E)$ associated to a hereditary subset $H$ of $E^0$ depends only on the saturation of $H$. This problem is solved in \cite[Theorem~5.3]{aHLRS2}, which gives a recipe for finding all the KMS$_\beta$ states of $\TT C^*(E)$ and $C^*(E)$ for fixed $\beta$.

\section{$C^*$-algebras from local homeomorphisms}\label{sec:local}

We consider a compact Hausdorff space $Z$ and a surjective local homeomorphism $h:Z\to Z$. In our main examples in the next section, $Z$ will be the infinite-path space $E^\infty$ of a finite directed graph with the topology inherited from the product space $(E^1)^\infty$, and $h$ will be the backward shift $\sigma$ defined by
\[
\sigma(e_1e_2e_3\cdots)=e_2e_3\cdots.
\]
If $E$ has no sources, then $\sigma$ is a homeomorphism on each cylinder set $Z(\mu)$, and hence is a local homeomorphism; if $E$ has no sinks, then $\sigma$ is also surjective. 
So we shall suppose in the rest of this paper that $E$ is a finite graph with no sinks or sources, and then $\sigma:E^\infty\to E^\infty$ is a good example to bear in mind for this section.

We can view $C(Z)$ as a Hilbert bimodule $X$ over the $C^*$-algebra $C(Z)$, by setting $(a\cdot x\cdot b)(z)=a(z)x(z)b(h(z))$ and 
\[
\langle x,y \rangle(z)=\sum_{h(w)=z}\overline{x(w)}y(w)\quad \text{for $x,y\in X$.}
\]
This Hilbert bimodule has both a Toeplitz algebra $\TT(X)$ and a Cuntz-Pimsner algebra $\OO(X)$: the Toeplitz algebra is generated by a representation $(\psi,\pi)$ characterised by $\psi(a\cdot x\cdot b)=\pi(a)\psi(x)\pi(b)$ and $\pi(\langle x,y \rangle)=\psi(x)^*\psi(y)$, and the Cuntz-Pimsner algebra \cite{P} is a quotient of $\TT(X)$. For our purposes, all the necessary background material is in Chapter~8 of \cite{R}. The Cuntz-Pimsner algebra  $\OO(X)$ is an example of Katsura's topological-graph algebras: in the conventions of \cite[Chapter~9]{R} (which are a little different from those in Katsura's original paper \cite{K}), the graph is $(Z,Z,\id,h)$.

The Toeplitz algebra $\TT(X)$ carries a gauge action $\gamma$ of the circle characterised by $\gamma_z(\psi(x))=z\psi(x)$ and $\gamma_z(\pi(a))=\pi(a)$, and this lifts to a dynamics $\alpha:t\mapsto \gamma_{e^{it}}$. The kernel of the quotient map onto $\OO(X)$ is invariant under $\gamma$, and hence we also get a dynamics  on $\OO(X)$ (still denoted by $\alpha$). 

Thomsen \cite{Th1} has studied the KMS states of the quotient system $(\OO(X),\alpha)$ (and he worked with much more general systems $(Z,h)$\,). He showed that the possible inverse temperatures of the KMS states all lie in a finite interval $[\beta_l,\beta_c]$, and gave formulas for upper and lower bounds:
\begin{align*}\label{defb}
\beta_c&=\limsup_{n\to \infty}\Big(n^{-1}\ln\Big(\max_{z\in Z}|h^{-n}(z)|\Big)\Big),\ \text{and}\\
\beta_l&=\limsup_{n\to \infty}\Big(n^{-1}\ln\Big(\min_{z\in Z}|h^{-n}(z)|\Big)\Big)
\end{align*}
(applying \cite[Theorem~6.8]{Th1} with the function $F\equiv 1$; see \cite[Remark~6.3]{AaHR} for the connections with Thomsen's notation). We are not aware that Thomsen has discussed the extent to which these bounds might be sharp.

In our recent work with Afsar \cite{AaHR}, we have studied the KMS states of the Toeplitz system $(\TT(X),\alpha)$. We viewed $C(Z)$ as a continuous analogue of the (finite-dimensional) space $C(E^0)$, and followed  the strategy of \cite{aHLRS1}. We found that, for inverse temperatures $\beta$ larger than Thomsen's $\beta_c$, the KMS$_\beta$ states are parametrised by a simplex $\Sigma_\beta$ of finite measures $\epsilon$ on $Z$ satisfying a normalisation condition of the form 
\[
\int f_\beta\,d\epsilon=1,
\]
where $f_\beta$ is a fixed continuous function defined by summing a series like that defining $y_\beta$ in \eqref{ybeta} \cite[Theorem~5.1]{AaHR}. At $\beta_c$, there is a phase transition: we can see by passing to limits as $\beta$ decreases to $\beta_c+$ that there exist KMS$_{\beta_c}$ states on $\TT(X)$, and can argue by mimicking our earlier  results in \cite{aHLRS1} that at least one of them factors through $\OO(X)$. (This is Theorem~6.1 in \cite{AaHR}.) So, in our generality at least, Thomsen's upper bound is sharp.

If $E$ is a finite graph, then there is a natural Hilbert bimodule $X(E)$ over the commutative $C^*$-algebra $C(E^0)$, and the Toeplitz algebra $\TT(X(E))$ was the original model of the Toeplitz-Cuntz-Krieger algebra $\TT C^*(E)$ (see \cite{FR} and \cite[Chapter~8]{R}). This bimodule is not given by a local homeomorphism, so it does not quite fit the set-up of the present section, but the analysis of \cite{AaHR} was  inspired by analogy with that of \cite{aHLRS1}. As we mentioned earlier, we can also directly apply the results of \cite{AaHR} to the shift $\sigma$ on the compact path space $E^\infty$, and this gives another connection to the results of \cite{aHLRS1} and \cite{aHLRS2}.

\section{Shifts on path spaces}\label{sec:shift}

We consider again a finite directed graph $E$ with no sinks or sources, and the infinite-path space $E^\infty$. Then $E^\infty$ is a compact Hausdorff space and the backward shift $\sigma:E^\infty\to E^\infty$ is a surjective local homeomorphism. So as in \S\ref{sec:local}, we can consider the Hilbert bimodule over the commutative $C^*$-algebra $C(E^\infty)$ with underlying space $X=C(E^\infty)$. At this point we choose to write $X(E^\infty)$ for $X$ to emphasise that this is not the graph bimodule $X(E)$ studied in \cite{FR} and \cite[Chapter~8]{R}.

The topology on $E^\infty$ arises from viewing it as a subset of the infinite product $(E^1)^\infty$ of the finite set $E^1$, and the cylinder sets
\[
Z(\mu)=\{x\in E^\infty: x_i=\mu_i\text{ for $i\leq |\mu|$}\} 
\]
associated to finite paths $\mu\in E^*$ form a basis of compact-open sets for the topology on $E^\infty$. Then a straightforward calculation shows:

\begin{lem}\label{includeTalgs}
The elements $P_v:=\pi(\chi_{Z(v)})$ and $S_e:=\psi(\chi_{Z(e)})$ of $\TT(X(E^\infty))$ form a Toeplitz-Cuntz-Krieger $E$-family.
\end{lem}

The universal property of the Toeplitz algebra $\TT C^*(E)$ now gives a homomorphism $\pi_{P,S}:\TT C^*(E)\to \TT (X(E^\infty))$. Corollary~4.2 of \cite{FR} implies that this homomorphism is injective, and it is equivariant for the gauge actions, and hence for the various dynamics $\alpha$ studied in \S\ref{sec:1} and \S\ref{sec:local}. So composing with $\pi_{P,S}$ takes KMS$_\beta$ states of $(\TT(X(E^\infty)),\alpha)$ to KMS$_\beta$ states of $(\TT C^*(E),\alpha)$. Now we have KMS$_\beta$ states of $\TT(X(E^\infty))$ for $\beta$ larger than Thomsen's $\beta_c$, and KMS$_\beta$ states of $\TT C^*(E)$ for $\beta>\ln \rho(A)$, where $A$ is the vertex matrix of $E$. We reconcile this in the following reassuring lemma, which is Proposition~7.3 of \cite{AaHR}.  (Note that the condition on $E$ is there to ensure that $\rho(A)>0$, so that $\ln \rho(A)$ makes sense.)

\begin{lem}\label{reassure}
Suppose that $E$ is a directed graph with at least one cycle. Then
\[
\frac{1}{n}\ln\Big(\max_{x\in E^\infty}|\sigma^{-n}(x)|\Big)\to \ln \rho(A)\quad\text{as $n\to \infty$.}
\]
\end{lem}

Thus Thomsen's $\beta_c$ is our $\ln\rho(A)$, and the range of possible $\beta$ in Theorem~\ref{paraKMS} is the the same as that in \cite[Theorem~5.1]{AaHR}. Suppose that $\beta>\ln\rho(A)$, that $\mu$ is a measure on $E^\infty$ satisfying the hypothesis $\int f_\beta\,d\mu=1$ of \cite[Theorem~5.1]{AaHR}, and that $\phi^\mu$ is the corresponding state of $\TT(X(E^\infty))$. Then Proposition~7.4 of \cite{AaHR} says that $\phi^\mu\circ\pi_{P,S}$ is the state $\phi_{\epsilon}$ of \cite[Theorem~2.1]{aHLRS1} associated to the vector $\epsilon=\epsilon(\mu)=\big(\mu(Z(v))\big)$ in $[0,\infty)^{E^0}$.

Every state $\phi_\epsilon$ of $(\TT C^*(E),\alpha)$ has the form $\phi^\mu\circ\pi_{P,S}$ for some measure $\mu$ on $E^*$ satisfying $\int f_\beta\,d\mu=1$ \cite[Corollary~7.6]{AaHR}. In the proof of this result, such a measure $\mu$ is constructed as a measure on the inverse limit $E^\infty=\varprojlim_nE^n$, and an examination of the construction shows that there is considerable leeway in building such a measure. Indeed, for each $\epsilon$ satisfying the normalisation relation $y^\beta\cdot\epsilon=1$ of \cite{aHLRS1},
\[
\Big\{\lambda\in M(E^\infty)_+:\int f_\beta\,d\lambda=1\text{ and }\lambda(Z(\lambda))=\epsilon_v \text{ for }v\in E^0\Big\}
\]
is a simplex of codimension $|E^0|+1$ in the cone $M(E^\infty)_+$ of positive measures. Thus there are many more KMS$_\beta$ states on $\TT(X(E^\infty))$ than on $\TT C^*(E)$. 

The injection $\pi_{P,S}:\TT C^*(E)\to \TT(X(E^\infty))$ is certainly not surjective --- if for no other reason, because $\TT(X(E^\infty))$ has many more KMS states.
However, Proposition~7.1 of \cite{AaHR} says that $\pi_{P,S}$ induces an isomorphism of the Cuntz-Krieger algebra $C^*(E)$ \emph{onto} $\OO(X(E^\infty))$! (This observation is essentially due to Exel \cite{E} and Brownlowe \cite{BRV}.) Since this isomorphism also intertwines the dynamics of \cite{aHLRS1} and that of \cite{AaHR}, the latter algebra has effectively the same KMS states as $C^*(E)$.

We now return to the dumbbell graph $E$
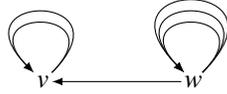
\[
\begin{tikzpicture}
    \node[inner sep=1pt] (v) at (0,0) {$v$};
    \node[inner sep=1pt] (w) at (2,0) {$w$};
    \draw[-latex] (w)--(v);
    \foreach \x in {0,2} {
        \draw[-latex] (v) .. controls +(1.\x,1.\x) and +(-1.\x,1.\x) .. (v);
    }
    \foreach \x in {0,2,4} {
        \draw[-latex] (w) .. controls +(1.\x,1.\x) and +(-1.\x,1.\x) .. (w);
    }
\end{tikzpicture}
\]
which we discussed in Example~\ref{2nddumb} of \S\ref{sec:dumb}. The system $(C^*(E),\alpha)$ has KMS$_\beta$ states for $\beta=\ln 3=\ln\rho(A)$ and $\beta=\ln 2=\ln\rho(A_{\{v\}})$. Thus so does $\OO(X(E^\infty))$. We have already seen in Lemma~\ref{reassure} that $\beta_c=\ln\rho(A)$ in general. For this $E$ and $x\in E^\infty$, we can compute
\[
|\sigma^{-n}(x)|=|E^nr(x)|=\begin{cases}
2^n&\text{if $r(x)=v$}\\
3^n+\sum_{j=0}^{n-1}3^j2^{n-1-j}&\text{if $r(x)=w$.}
\end{cases}
\]
Thus $\min_x|\sigma^{-n}(x)|=2^n$ is attained when $r(x)=v$, and Thomsen's $\beta_l$ is $\ln 2$. So for the  local homeomorphism $\sigma:E^\infty\to E^\infty$, the lower bound $\beta_l$ in \cite[Theorem~6.8]{Th1} is also sharp.

It is easy to see with dumbbell graphs that there can be KMS$_\beta$ states at inverse temperatures strictly between $\beta_l$ and $\beta_c$. For example, with $E$ the following graph
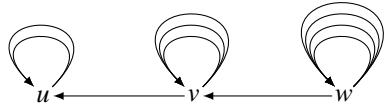
\[
\begin{tikzpicture}
    \node[inner sep=1pt] (u) at (-2,0) {$u$};
    \node[inner sep=1pt] (v) at (0,0) {$v$};
    \node[inner sep=1pt] (w) at (2,0) {$w$};
    \draw[-latex] (w)--(v);
    \draw[-latex] (v)--(u);
    \foreach \x in {0,2} {
        \draw[-latex] (u) .. controls +(1.\x,1.\x) and +(-1.\x,1.\x) .. (u);
    }
    \foreach \x in {0,2,4} {
        \draw[-latex] (v) .. controls +(1.\x,1.\x) and +(-1.\x,1.\x) .. (v);
    }
     \foreach \x in {0,2,4,6} {
        \draw[-latex] (w) .. controls +(1.\x,1.\x) and +(-1.\x,1.\x) .. (w);
    }
\end{tikzpicture}
\]
both $C^*(E)$ and $\OO(X(E^\infty))$ have KMS states at inverse temperatures $\ln2$, $\ln 3$ and $\ln 4$.

There are, however, interesting constraints on the possible inverse temperatures $\beta$. First, since $e^\beta$ has to be the spectral radius of an irreducible integer matrix, it has to be an algebraic number. But there are also other, more subtle constraints. The issue is discussed, along with relevant results of Lind \cite{L}, in \cite[\S7.1]{aHLRS2}.

\begin{acknowledgement}
This research was supported by the Marsden Fund of the Royal Society of New Zealand. We also thank our collaborators Zahra Afsar, Marcelo Laca and Aidan Sims, the referee for some constructive suggestions, and the organisers of the Abel symposium for a marvellous experience.\end{acknowledgement}

\end{document}